%% file: Meromorphic_Ext.tex
\documentclass[11,letterpaper,leqno]{article}
\usepackage[a4paper, hmargin=2.5cm, vmargin={2.5cm, 2.5cm}]{geometry}
\usepackage{enumitem}
\usepackage{amsfonts}
\usepackage{amssymb}
\usepackage{amsmath}
\usepackage{amsthm}
\usepackage{graphicx}
\usepackage[hyperindex]{hyperref}
\usepackage{mathrsfs}
\usepackage{subfiles}
\usepackage{xcolor}

\def\R{\ensuremath{\mathbb R} }
\def\C{\ensuremath{\mathbb C} }
\def\Q{\ensuremath{\mathbb Q} }

\def\N{\ensuremath{\mathbb N} }
\def\Z{\ensuremath{\mathbb Z} }

\def\eps{\ensuremath{\varepsilon} }
\def\vect{\ensuremath{\overset{\rightharpoonup}}}
\newcommand \nin{\not\in}

\newcommand* \bb [1]{\left({#1}\right)}
\newcommand* \sbb [1]{\left[{#1}\right]} 
\newcommand* \bset [1]{\left\{{#1}\right\}} 
\newcommand* \limit [2]{\underset{{#1}\rightarrow{#2}}{\lim}\;}
\newcommand* \bunion[3]{\bigcup_{{#1} = {#2}}^{#3}}

\newcommand* \sumit [3]{\underset{{#1} = {#2}}{\overset{#3}\sum}\;}
\newcommand* \prodit [3]{\underset{{#1} = {#2}}{\overset{#3}\prod}\;}
\newcommand* \abs[1] {\left|{#1}\right|}

\theoremstyle{definition}
\newtheorem{thm}{Theorem}[section]
\newtheorem*{thm*}{Theorem}

\newtheorem{lem}[thm]{Lemma}

\newtheorem{defn}[thm]{Definition}
\newtheorem{rmk}[thm]{Remark}
\newtheorem{obs}[thm]{Observation}

\makeatletter
\newcommand\RedeclareMathOperator{%
  \@ifstar{\def\rmo@s{m}\rmo@redeclare}{\def\rmo@s{o}\rmo@redeclare}%
}
\newcommand\rmo@redeclare[2]{%
  \begingroup \escapechar\m@ne\xdef\@gtempa{{\string#1}}\endgroup
  \expandafter\@ifundefined\@gtempa
     {\@latex@error{\noexpand#1undefined}\@ehc}%
     \relax
  \expandafter\rmo@declmathop\rmo@s{#1}{#2}}
\newcommand\rmo@declmathop[3]{%
  \DeclareRobustCommand{#2}{\qopname\newmcodes@#1{#3}}%
}
\@onlypreamble\RedeclareMathOperator
\makeatother

\RedeclareMathOperator{\Re}{Re} 

\AtEndDocument{\bigskip{\footnotesize%
  \addvspace{\medskipamount}
  (A.~Gl\"ucksam) \textsc{Einstein Institute of Mathematics, The Hebrew University of Jerusalem, Jerusalem ,Israel} \par
  \textit{E-mail address}: \texttt{adi.glucksam@math.huji.ac.il}
    \addvspace{\medskipamount}
  (Y.J.~Zou) \textsc{Department of Mathematics and Statistics, Oakland University, Rochester, MI, USA} \par
  \textit{E-mail address}: \texttt{yzou@oakland.edu}

}}

\begin{document}
\title{A Comparison Test for Meromorphic Extensions}

\author{Adi Gl\"ucksam and Yuzhou Joey Zou}

\maketitle

\begin{abstract}
We provide a comparison test for meromorphic extensions, i.e., if two series are  ``close enough" then the existence of a meromorphic extension of one to the entire complex plane ensures a similar extension for the other. We use this result to generate new examples of Dirichlet series admitting meromorphic extensions. Moreover, we demonstrate that our requirements are optimal by constructing a collection of counterexamples where the series are ``close but not enough": one series admits a meromorphic extension while the other possesses a natural boundary.

\end{abstract}

\section{Introduction}
\subfile{sections/introduction}

\section{Decomposition of Series}
\subfile{sections/pre}

\section{The proof of Theorem \ref{thm:extension}}
\subfile{sections/proof}

\section{Examples and Counterexamples}
\subfile{sections/examples}

\end{document}

%% file: sections/introduction.tex
Let $\bset{c_n}\subset(0,\infty)$ be a sequence satisfying $\limit n\infty c_n=0$ and let $\bset{a_n}\subset\C$. Define the \emph{Dirichlet series}
$$
D(s):=\sumit n 1\infty a_n\cdot c_n^s.
$$
Assume that there exists $S_0$ so that the Dirichlet series above converges absolutely in the half-plane $\mathcal{P}_0:=\bset{s\in\C\,:\,\Re(s)>S_0}$ (the smallest such $S_0$ is sometimes called the \emph{abscissa of absolute convergence} for the series).

We say the series has a \emph{meromorphic extension to the entire complex plane} if there exists $g\in Mer(\C)$ so that $g|_{\mathcal{P}_0}\equiv \sumit n 1\infty a_n\cdot c_n^s$. If this is not the case it could be, for example, that the function has an essential singularity or that the singular points of the function are not isolated. In the latter case we say it has a \emph{natural boundary}. Formally
\begin{defn}
Let $U\subset\C$ be an open set, and let $f:U\rightarrow\C$ be meromorphic in $U$. We say $f$ has a \emph{natural boundary} along a non-trivial arc, $A\subset\partial U$, if for every $z\in A$,
there is no function, $g$, meromorphic in some neighbourhood of $z$, $N_z$, satisfying $g(w)=f(w)$ for all $w\in N_z\cap U$.
\end{defn}

A fundamental example where a Dirichlet series admits a meromorphic extension is when $a_n=1$ and $c_n=a^{-n}$. In this case, the Dirichlet series is a geometric progression, and the sum $\sumit n 1 \infty (a^{-n})^s$ converges on $\mathcal{P}_0 = \{s\in\C\,:\, \Re(s)>0\}$ to $\frac{1}{a^s-1}$ if $a>1$. The function $s\mapsto\frac1{a^s-1}$ is holomorphic in $\C\setminus\bset{ \frac{2\pi ik}{\log(a)},\; k\in\Z}$ and has simple poles. In particular, it is 
holomorphic on $\mathcal{P}_0$, where it identifies with the original Dirichlet sum.

On the other hand, not all Dirichlet series admit meromorphic extensions. For example, by Hadamard's gap theorem (see e.g.\ \cite{R}, Section 11.2), the series $\sumit j 0 \infty z^{p_j}$ has a natural boundary at $\bset{\abs z=1}$ if $p_j$ is a sequence of positive integers such that $\inf_j \frac{p_{j+1}}{p_j} > 1$. Using the notation above, defining $a_n = 1$ if $n = p_j$ for some $j$, and $a_n=0$ otherwise, we see that the series $\sumit n 0 \infty a_n\cdot z^n$ has a natural boundary at $\bset{|z|=1}$; substituting $z = e^{-s}$shows that the series $\sumit n 0 \infty a_n\cdot e^{-ns}$ has a natural boundary at $\bset{s\in\C\,:\,\Re(s)=0}$.

We thus ask the following question: assume that $\sumit n 1 \infty a_n\cdot c_n^s$ has a meromorphic extension to the entire complex plane, and let $\widehat{c}_n$ be some perturbation of $c_n$. When does $\sumit n 1 \infty a_n\cdot\widehat{c}_n^{\;s}$ have a meromorphic extension to the entire complex plane as well?

This question asks whether there exists a form of a comparison test for meromorphic extensions. The existence of a meromorphic extension of a given series is generally difficult to prove and relies on some specific structures and symmetries given by the series. It is therefore helpful to know, given a series known to admit a meromorphic extension, when one can leverage this fact to give extensions for ``similar-looking'' series, i.e., perturbations of the known series. Recall that the classical comparison test tells us that if two sequences, $\bset{a_n},\bset{b_n}$, satisfy $\frac1C\le \frac{\abs{a_n}}{\abs{b_n}}\le C$ for some $C>1$ and for all $n$, then their corresponding series converge and diverge (in absolute value) together.  It allows us to reduce the problem of proving absolute convergence of a complex series by comparing it to a simpler, well-understood series. We aim to create a similar framework for meromorphic extensions. A notable example to keep in mind is the geometric series we saw earlier. Our results will allow one to leverage the existence of meromorphic extension for geometric series to obtain meromorphic extensions to series such as $\sumit n 1 \infty\bb{2^{-n}+3^{-n}}^s$ (see Example \ref{thm:examp_23}). Naturally, the conditions required to compare the existence of meromorphic extensions are stronger than those required to compare simple convergence (see section \ref{sec:res} below). Nevertheless, we show these requirements are necessary (see Theorem \ref{thm:example}).

Note that the ``classical'' setting with increasing exponents, i.e., $\sumit n 1 \infty a_n\cdot\lambda_n^{-s}$ with $\lambda_n\nearrow\infty$ can be represented in the setup described above by setting $c_n=\frac1{\lambda_n}$.

\subsection{Some History}
We outline several results that are related to the question presented above.
\begin{itemize}[label=$\bullet$]
\item Mandelbrojt \cite{M} established several foundational results relating the arithmetic properties of the exponents in a Dirichlet series to its analytic behavior, particularly concerning the location of singularities and the nature of its analytic continuation. Among these are ``comparison'' results for Hadamard-type compositions, which involve series of the form $\sumit n 1 \infty a_n\cdot b_n \cdot e^{-\lambda_n \cdot s}$ derived from two series with identical exponents. These coefficient-based modification results are similar to the questions we investigate, although in our case we focus on perturbing the bases $c_n$ (corresponding to the exponents $\lambda_n$ in Mandelbrojt's notation) rather than the coefficients.
\item Navas, Ruiz, and Varona \cite{NRV} 
related the convergence of a power series $\sumit n 0 \infty a_n\cdot z^n$ with the extendibility of the Dirichlet series $\sumit n 1 \infty a_{n-1}\cdot n^{-s}$, by essentially using the Mellin transform. Note that here $c_n = \frac{1}{n}$, and the coefficients, $a_n$, are the same coefficients as in the original power series. They also showed that changing the coefficients $a_n$ in a specific way preserved the property of having a meromorphic extension (as before, with $c_n=\frac1n$). For more information see Corollary 1 and 2 in \cite{NRV}.
\item Bhowmik and Schlage-Puchta showed in \cite{BSP} that if the Riemann zeta function has infinitely many zeroes outside the line $\frac12+it$, then a certain product of dilated Riemann zeta function has a natural boundary along the line $\{\Re(s)=\frac12\}$,  by showing that zeros accumulate (and are not cancelled by poles) as one approaches the purported natural boundary. They used this to show that certain Euler products (which can be written in terms of Dirichlet series) have natural boundaries. In addition, they proved a corresponding result regarding a randomized version of the modification.
\item Breuer and Simon investigated natural boundaries of power series in \cite{BS}, using the idea of \emph{right limits} from the spectral theory of Jacobi matrices. They used the result to show that a \emph{generic} (in either a probabilistic or topological sense) power series with radius of convergence $1$ will have a natural boundary at the unit circle.

\item Ram Murty and Sinha \cite{RS} studied meromorphic continuation of \emph{multiple Hurwitz zeta functions}, which are multivariable generalizations of the Hurwitz zeta function, initially defined on a subset of $\mathbb{C}^d$, to all of $\mathbb{C}^d$. A key technique in their paper, when specialized to one variable, recovers the meromorphic extension of the Hurwitz zeta function $\sum_{n=0}^{\infty}(n+a)^{-s}$ by expanding the summand in terms of $n^{-s}$ by using the binomial formula and invoking known properties of the Riemann zeta function.  This can be viewed as a special case of our Theorem \ref{thm:extension} and a special case of Theorem \ref{thm:zeta} (see also Remark \ref{rmk:hurwitz}).
 This technique was generalized in a recent paper by Wang, Mehta, and Kanemitsu \cite{WMK}, who studied perturbations of a generalized Dirichlet series $\sumit n 1 \infty a_n\cdot\lambda_n^{-s}$ to $\sumit n 1 \infty a_n\cdot(\lambda_n+a)^{-s}$ by applying similar methods.
\item 
Finally, the ideas used in the proof of Theorem \ref{thm:extension} below, showing meromorphicity of a perturbed Dirichlet series, are similar to the 
\emph{Factorization Method} used to establish meromorphic extension of Euler products. In this method, one factors out certain terms which combine to form known functions (namely, the Riemann zeta function), whose analytic properties are known, leaving a remainder which is analytic on a set larger than initially expected, and, in turn, establishes meromorphicity on a set larger than initially expected. See \cite{Ahmetaj} and \cite{Alberts} for expository sources on this method.
\end{itemize}

While the list above is quite extensive, these works primarily focus on Dirichlet series $\sum_{n=1}^{\infty} a_n \cdot c_n^s$ where the base, $c_n$, is either a fixed well-studied sequence (such as $\frac{1}{n}$ or $e^{-n}$) or a simple modification thereof (such as shifting by a constant, as in \cite{RS, WMK}), and they study how different coefficients, $a_n$, or specific structured perturbations of the base determine the existence---or absence---of a meromorphic extension. We are interested in a more general question: what happens if we perturb the \emph{base} $c_n$ in a controlled but otherwise arbitrary way, while leaving the coefficients $a_n$ unchanged. This is a problem that, to the best of our knowledge, has not been explored.

\subsection{A Comparison Test}\label{sec:res}
\subsubsection{On the Type of Perturbations}
In considering a perturbation, $\sumit n 1 \infty a_n\cdot\widehat{c}_n^{\;s}$, of a series, $\sumit n 1 \infty a_n\cdot c_n^s$, we will consider the case where $\widehat{c}_n$ can be written as a (nice) function of $c_n$. That is, we will assume
\[\widehat{c}_n = g(c_n),\]
where $g(z)$ is a function which admits an asymptotic expansion
\[g(z)\sim \alpha_0z^{\rho_0}+\alpha_1z^{\rho_1}+\dots,\]
where $\alpha_i\in\C$, $\alpha_0>0$ and $0<\rho_0<\rho_1<\dots$ is an increasing sequence. Without loss of generality, we may assume $\alpha_0=1$, since 
we can factor out $\alpha_0^s$ as an entire factor. Moreover, it will be convenient to factor out the leading order power $z^{\rho_0}$ out of $g$. Thus, if we set
\[\sigma = \rho_0,\quad \sigma_j = \rho_j-\rho_0>0\quad \text{ for }j\ge 1,\]
 we can write
 \[g(z) = z^\sigma(1+h(z)),\quad\text{ where }\quad h(z)\sim\alpha_1z^{\sigma_1}+\alpha_2z^{\sigma_2}+\dots\]
See the hypotheses of Theorem \ref{thm:extension} below for a precise definition of $h$ satisfying this asymptotic expansion. In this case, $\widehat{c}_n = c_n^\sigma(1+h(c_n))$, and the perturbed series becomes
\[\sumit n 1 \infty a_n\cdot\widehat{c}_n^{\;s} = \sumit n 1 \infty a_n\cdot c_n^{\sigma\cdot s}(1+h(c_n))^s,\]
which converges absolutely in $\left\{s\in\C\,:\,\Re(s)>\frac{S_0}\sigma\right\}$.

\subsubsection{Results}
We prove two complementary results---a positive result and a counterexample. To phrase our results, we define a \emph{bounded sector} to be a set of the form 
$\Gamma = \{r\cdot e^{it}\,:\, r\in\sbb{0,\eps},t\in\sbb{t_1,t_2}\}$ for some numbers $\eps,t_1,t_2$, with $\eps>0$ and $0\le t_2-t_1<2\pi$. For example, $\Gamma$ could be a line segment containing the origin, or closed bounded sector of any angle less than $2\pi$.

\begin{thm}\label{thm:extension}
Let $D(s):=\sumit n 1 \infty a_n\cdot c_n^s$ be a Dirichlet series with positive bases $c_n$ that converges absolutely in the half-plane $\mathcal{P}_0=\bset{\Re(s)>S_0}$, and extends meromorphically to $\C$. 
Let $\Gamma$ be a bounded sector containing the set $\{c_n\}$, and let $h:\Gamma\rightarrow\C$
be a function satisfying that for every $N$, there exists $C_N>1$ such that 
\begin{align}\label{eq:cond_h}
\abs{h(z)-\sumit j1N \alpha_j\cdot  z^{\sigma_j}}\le C_N\cdot \abs z^{\sigma_{N+1}} \text{ for all }z\in \Gamma 
\end{align}
where $\{\alpha_j\}\subset\C$ and $\bset{\sigma_j}\subset(0,\infty)$ satisfies $\sigma_j\nearrow \infty$.

Then the sum $\sumit n 1 \infty a_n\cdot \left(c_n^{\sigma}\bb{1+h(c_n)}\right)^s$ extends meromorphically to $\C$ as well.
\end{thm}
In other words, if $h(z)$ admits an asymptotic expansion in terms of (not necessarily integer) powers of $z$, then the property of meromorphic extension for the original series applies to the perturbation as well. 
Note that the function $h$ is not required to be smooth or even differentiable in a neighbourhood of the origin (see, e.g., the example in Section \ref{subsec:rough}).

Moreover, on one hand, condition \eqref{eq:cond_h} implies that $h(z)$ must be bounded on $\Gamma$ (which is bounded by definition), and on the other hand, if $h(z)$ is bounded, then it is enough to verify condition \eqref{eq:cond_h} only for $\abs z<\eps_N$ for some $\eps_N>0$, as we could always replace $C_N$ by $\frac C{\eps_N^{\sigma_{N+1}}}$ to accommodate all $z\in\Gamma$, if needed.

\begin{rmk}
While the bases, $\bset{c_n}$, are assumed to be positive, their perturbed counterparts, $\bset{\widehat{c_n}}$, are not assumed to even be real. In fact, Theorem \ref{thm:extension} holds even if the bases, $\bset{c_n}$, are not positive, but complex numbers. In this case, since $\bset{c_n}$ is a countable set, there exists a branch of the logarithm such that $\log(c_n)$ is well defined for all the elements in $\bset{c_n}$ and the sums
$$
\sumit n 1 \infty a_n\cdot c_n^s=\sumit n 1 \infty a_n\cdot \exp\bb{s\cdot\log(c_n)}
$$
are well defined and converge for some $s\in\C$. We write the proof for $c_n>0$ as otherwise many of the details become cumbersome to write, but the idea of the proof is the same. In fact, we only need that $h$ is defined on the set $\bset{c_n}$. Nevertheless, we assume that $h$ is defined on some 
bounded sector
to make the transition to complex bases simpler.
\end{rmk}

\begin{rmk}
Note that these requirements allow a large collection of perturbations including perturbations that are finite sums of powers (in this case one can set $\sigma_{N+k}=\sigma_N+k, a_{n+k}=0$ for all $k\in\N$) and perturbations that are not at all of this form, but are small enough like $h(z)=\exp\bb{-\frac{1}{\abs z}}$, which is smooth on $\bset{\abs z>0}$ but is not differentiable in the complex sense.
\end{rmk}

We show that the \emph{asymptotic expansion} requirement (condition \eqref{eq:cond_h}) on $h$ is truly necessary, and that a polynomial bound is not enough, in the following counterexample:

\begin{thm}\label{thm:example}
For every positive integer, $m$, there exists a function, $h(z)$, satisfying:
\begin{enumerate}\item $h$ is smooth on $(0,\eps)$.
\item $|h(z)|\le |z|^m$ on $\sbb{0,\eps}$.
\item the Dirichlet series $\sumit n 1 \infty \left(e^{-n}(1 + h(e^{-n}))\right)^s$ converges in $\bset{\Re(s)>-m}$, but has a natural boundary along the vertical line $\bset{\Re(s) = -m}$.
\end{enumerate}
\end{thm}

\subsection{Structure of the Paper}
The first part of the paper shows how given a finite approximation for $h$ by a power series (with real powers), the sum $\sumit n 1 \infty a_n\cdot c_n^s\bb{1+h(c_n)}^s$ decomposes into two sums: a main term composed of weighted sums of translations of the original sum, and an error term which is holomorphic in some half-plane growing with the highest power in the finite approximation. The second section provides the proof of Theorem \ref{thm:extension} using the decomposition presented in the first section. The last section is composed of two subsections. The first explores the usefulness of Theorem \ref{thm:extension} and proves several generic cases where the theorem can be applied. In the second subsection we construct the counterexample showing the requirements in Theorem \ref{thm:extension} are necessary.

\subsection{Acknowledgements}
The authors would like to thank Nir Avni for motivating the question, and several useful discussions. The authors also thank M.\ Sodin, S.\ Sodin, and A.\ Eremenko for suggesting looking at natural boundaries.
The authors gratefully acknowledge support from the Department of Mathematics at Northwestern University, where a significant part of the discussion and work on this project was conducted.
The first author is grateful for the support of the Golda Meir fellowship.

%% file: sections/pre.tex
Given a finite approximation for $h$ by a power series (with real powers), we show that the sum \\ ${\sumit n 1 \infty a_n\cdot c_n^s\bb{1+h(c_n)}^s}$ decomposes into two sums: a main term composed of weighted sums of translations of the original sum, and an error term which is holomorphic in some half-plane growing with the highest power in the finite approximation.
\paragraph*{Notation}
Throughout the paper we shall use the following notation:
\begin{enumerate}
\item $(s)_k$ will denote the falling factorial defined by $(s)_k:=\prodit j 0 {k-1}\bb{s-j}$.
\item We denote by $\N_0$  the set of all non-negative integers, i.e., $\N_0:=\N\cup\bset 0$.
\item We denote by $B(0,\eps)$ the open disk of radius $\eps$ in $\C$ centered at the origin, namely, $B(0,\eps):=\bset{z\in\C\;\colon\;\abs z<\eps}$, and by $\Gamma$ the bounded sector introduced above Theorem \ref{thm:extension}.
\item We write $\vect k$ for a multi-index $\vect k=\bb{k_1,\cdots,k_{m}}\in\N^m$, and define the multinomial coefficients
$$
\bb{ \sumit j 1 m k_j\atop { k_1,k_2,\cdots,k_m}}=\frac{\bb{\sumit j 1 m k_j}!}{\prodit j 1 m k_j!}.
$$
\item Given a multi-index, $\vect k=(k_1,\cdots,k_m)$, and a sequence of real numbers, $\bset{\alpha_j}_{j=1}^N$, where $N\ge m$, we define
$$
A\bb{\vect k,\bset{a_j}}:=\bb{ \sumit j 1 m k_j\atop { k_1,k_2,\cdots,k_m}}\prodit j1 m\alpha_j^{k_j},\quad \text{and}\quad \vect{\sigma}\cdot\vect{k} := \sumit j 1 m\sigma_{j}\cdot k_j.
$$
\item With $\sigma$, $\{\sigma_j\}$ the parameters appearing in the hypotheses of Theorem \ref{thm:extension}, and $S_0$ the abscissa of absolute convergence for $D(s)$ as defined in the introduction, define the half-planes $P_N=\bset{s\in\C\,:\, \Re(s)> \frac{S_0-\sigma_{N+1}}{\sigma}}$. Note that
\[s\in P_N\iff \sigma\cdot s+\sigma_{N+1}\in\mathcal{P}_0,\]
where $\mathcal{P}_0=\bset{s\in\C, \Re(s)>S_0}$ is the half-plane where the original sum is absolutely convergent and converges to a holomorphic map. By convention, we set $\sigma_0=0$, and let $P_{-1} = \left\{s\in\C\,:\,\Re(s)>\frac{S_0}\sigma\right\}$; this is the half-plane where the perturbed series will converge absolutely.
\end{enumerate}
Our first lemma concerns expanding $(1+h(z))^s$ when $h$ admits a finite expansion in powers of $z$, up to an error:

\begin{lem}\label{lem:apprx_poly_h}
Let 
$\Gamma$ be a bounded sector, and $h:\Gamma\rightarrow\C$
be a function satisfying $h(z)=\sumit j1N \alpha_j\cdot  z^{\sigma_j}+E(z)$, where $\bset{\sigma_j}\subset(0,\infty)$ is a monotone increasing sequence and $\abs{E(z)}\le C\cdot\abs z^{\sigma_{N+1}}$ for some $C>0$. We define  
$$
\mathcal E_k(z):=\sumit {k_{N+1}}1 k\frac{\bb{E(z)}^{k_{N+1}}}{(k_{N+1})!}\underset{\vect k, \sumit j 1 {N} k_j=k-k_{N+1}}\sum A\bb{\vect k,\bset{\alpha_j}}\cdot z^{\vect{\sigma}\cdot\vect{k}}.
$$
Fix $\delta\in\bb{0,\frac12}$ such that
\begin{equation}
\label{eq:delta}
\sumit j 1 N \abs{\alpha_j}\cdot\delta^{\sigma_j}<\frac12\quad\quad\text{and}\quad\quad C\cdot \delta^{\sigma_{N+1}}<\frac12.
\end{equation}
Then, whenever $z\in\Gamma$ satisfies $\abs z<\delta$,
\begin{enumerate}[label=(\roman*)]
\item\label{eq:1+h-exp} 
$\bb{1+h(z)}^s=1+\sumit k 1\infty\frac{(s)_k}{k!}\bb{\underset{\vect k,\; \sumit j 1 {N} k_j=k}\sum A\bb{\vect k,\bset{\alpha_j}}z^{\vect{\sigma}\cdot\vect{k}}+\mathcal E_k(z)}.$
\item \label{item:err}
$\abs{\mathcal E_k(z)}\le C\cdot  \abs z^{\sigma_{N+1}}\cdot k\cdot 2^{-(k-1)}.$
\end{enumerate}
\end{lem}
\begin{rmk}
Note that $\delta$ satisfying \eqref{eq:delta} always exist as the requirement that $\sigma_j>0$ implies that the left hand side of both inequalities tend to 0 as $\delta\searrow0$.
\end{rmk}
\begin{proof}
Recall that the Taylor expansion of the function $z\mapsto(1+z)^s$ centred at $z=0$ is
$$
\bb{1+z}^s=1+\sumit k 1 \infty \frac{(s)_k}{k!}\cdot z^k.
$$
Note that if $\abs z<\frac12$, then we may use the principal branch of the logarithm to define $z\mapsto \bb{1+z}^s$.

We therefore aim to approximate $\bb{h(z)}^k$ by sums of powers. Following the multinomial theorem, for every $m,n\in\N$
$$
\bb{x_1+x_2+\cdots+x_m}^n=\underset{\vect k, \sumit j 1 m k_j=n}\sum {n\choose k_1,k_2,\cdots,k_m}\prodit j 1 m x_j^{k_j}=\underset{\vect k, \sumit j 1 m k_j=n}\sum A\bb{\vect k,\bset{x_j}},
$$
using our notation.

In particular, for every $k$
\begin{align*}
\bb{h(z)}^k=\bb{\sumit j 1{N}\alpha_j\cdot z^{\sigma_j}+E(z)}^k=\underset{\vect k,\; \sumit j 1 {N+1} k_j=k}\sum {k\choose k_1,k_2,\cdots,k_{N+1}}\prodit j 1 {N} \alpha_{j}^{k_j}\cdot z^{\sigma_{j}\cdot k_j}\cdot \bb{E(z)}^{k_{N+1}}\\
=\underset{\vect k,\; \sumit j 1 { N}k_j=k}\sum A\bb{\vect k,\bset{\alpha_j}}z^{\vect{\sigma}\cdot\vect{k}}+\mathcal E_k(z).
\end{align*}
Combining this with the Taylor expansion of $\bb{1+z}^s$, we see that
$$
\bb{1+h(z)}^s=1+\sumit k 1\infty\frac{(s)_k}{k!}\cdot \bb{h(z)}^k=1+\sumit k 1\infty\frac{(s)_k}{k!}\bb{\underset{\vect k,\; \sumit j 1 {N} k_j=k}\sum A\bb{\vect k,\bset{\alpha_j}}z^{\vect{\sigma}\cdot\vect{k}}+\mathcal E_k(z)},
$$
concluding the proof of part \ref{eq:1+h-exp}.

To see the proof of part \ref{item:err} first note that using the multinomial theorem again
\begin{equation}\label{eq:multi}
\begin{aligned}
\underset{\vect k, \sumit j 1 {N} k_j=k-k_{N+1}}\sum\abs{A\bb{\vect k,\bset{\alpha_j}}z^{\vect{\sigma}\cdot\vect{k}}}=
\underset{\vect k, \sumit j 1 {N} k_j=k-k_{N+1}}\sum A\bb{\vect k,\bset{|\alpha_j|}}\cdot \abs z^{\vect{\sigma}\cdot\vect{k}}\quad\quad\quad\quad\quad\quad\quad\quad\quad\\
= \bb{\abs z^{\sigma_1}\cdot\abs{\alpha_1}+\cdots+\abs z^{\sigma_N}\cdot\abs{\alpha_N}}^{k-k_{N+1}}\le \bb{\sumit j 1 N \abs{\alpha_j}\delta^{\sigma_j}}^{k-k_{N+1}}\le 2^{-(k-k_{N+1})}.
\end{aligned}
\end{equation}
This implies that since $\abs{E(z)}\le C\cdot \abs z^{\sigma_{N+1}}\le C\cdot \delta^{\sigma_{N+1}} \le 1/2$,
\begin{align*}
\abs{\mathcal E_k(z)}&\le\sumit {k_{N+1}}1 k\abs{\bb{E(z)}^{k_{N+1}} \underset{\vect k, \sumit j 1 {N} k_j=k-k_{N+1}}\sum A\bb{\vect k,\bset{\alpha_j}}\cdot z^{\vect{\sigma}\cdot\vect{k}}}\\
&\le C|z|^{\sigma_{N+1}}\sumit {k_{N+1}} 1 k (C|z|^{\sigma_{N+1}})^{k_{N+1}-1}\cdot 2^{-(k-k_{N+1})} \\
&\le C|z|^{\sigma_{N+1}}\sumit {k_{N+1}} 1 k 2^{-(k_{N+1}-1)}\cdot 2^{-(k-k_{N+1})} = C|z|^{\sigma_{N+1}}\cdot k\cdot 2^{-(k-1)},
\end{align*} 
concluding the proof of the lemma.
\end{proof}

We now show that the portion of the Dirichlet series arising from the ``error'' term in expanding $(1+h(z))^s$ is holomorphic on a large half-plane:
\begin{lem}\label{lem:analytic_tail}
Let $\Gamma$ be a bounded sector containing the set $\bset{c_n}$, and let $h:\Gamma\to\C$ be a function of the form $h(z) = \sumit j 1 N \alpha_j\cdot z^{\sigma_j} + E(z)$, where $\{\sigma_j\}$, $E(z)$ satisfy the requirements of Lemma \ref{lem:apprx_poly_h}.
Let $\delta$ 
satisfy the conditions in \eqref{eq:delta}, and suppose that the sum $D(s):=\sumit n 1 \infty a_n\cdot c_n^s$ converges to a holomorphic function in $\mathcal{P}_0=\bset{\Re(s)>S_0}$.
For every $\nu$ satisfying that
$$
\abs{c_n}<\delta\quad\text{ for all }n\ge \nu
$$
the sum 
$$
\sumit k 1\infty\frac{(s)_k}{k!}\sumit n {\nu} \infty a_n\cdot c_n^{\sigma\cdot s}\cdot\mathcal E_k(c_n),
$$
is absolutely convergent in $P_N$ and converges to a function which is holomorphic in $P_N$, where the terms $\mathcal E_k(z)$ are as described in Lemma \ref{lem:apprx_poly_h} (note that these depend on $h$).
\end{lem}
\begin{rmk}
Condition \eqref{eq:delta} changes with the constant $N$. If the asymptotic expansion of $h$ has infinitely many terms, then the smallest $\nu$ that can be chosen in the lemma above will also depend on the chosen finite approximation order, $N$.
\end{rmk}
\begin{proof}
We would like to bound the sum 
$$
\sumit k 1\infty\abs{\frac{(s)_k}{k!}}\sumit n {\nu} \infty \abs{a_n\cdot c_n^{\sigma\cdot s}\cdot\mathcal E_k(c_n)}
$$
locally uniformly in $P_N$.

We note that since $1+x\le e^x$ for all real $x$,
\begin{equation}
\begin{aligned}
\label{eq:falling}
\abs{\frac{(s)_k}{k!}} = |s|\prod_{j=1}^{k-1}\frac{|s-j|}{j+1} &\le |s|\prod_{j=1}^{k-1}\left(1+\frac{|s+1|}{j+1}\right) \\
&\le |s|\exp\left(\sumit j 1 {k-1} \frac{|s+1|}{j+1}\right)\le |s|\exp(|s+1|\log(k)) = |s|k^{|s+1|}.
\end{aligned}
\end{equation}
While the manipulation requires $k\ge 2$ to make sense, the result is true for $k=1$ as well. Combining this with part \ref{item:err} of Lemma \ref{lem:apprx_poly_h}, we obtain
\begin{align*}
\sumit n {\nu} \infty\sumit k 1\infty\abs{\frac{(s)_k}{k!}} \abs{a_n\cdot c_n^{\sigma\cdot s}\cdot\mathcal E_k(c_n)}&\le\abs s\sumit n {\nu} \infty\sumit k 1\infty k^{\abs {s+1}} \abs{a_n}\cdot \abs{c_n^{\sigma\cdot s}}\cdot  C\cdot \abs {c_n}^{\sigma_{N+1}}\cdot k\cdot 2^{-(k-1)}\\
&\le C\cdot\abs s\sumit n {\nu} \infty\abs{a_n}\cdot \abs{c_n^{\sigma\cdot s+\sigma_{N+1}}}\sumit k 1\infty k^{|s+1|+1} \cdot  2^{-k+1}<\infty
\end{align*}
whenever $\sigma\cdot s+\sigma_{N+1}\in \mathcal{P}_0$, i.e.,\ when $s\in P_N$, concluding the proof of the lemma.
\end{proof}

Finally, we establish the meromorphicity of the ``main term'':
\begin{lem}\label{lem:mero-comb}
Under the assumptions of Lemma \ref{lem:analytic_tail}, let $D_\nu(s):=\sumit n {\nu} \infty a_n\cdot c_n^s$. 
Then the series
$$
\widetilde D(s):=\sumit k {0}\infty\frac{(s)_k}{k!}\underset{\vect k,\; \sumit j 1 {N} k_j= k}\sum A\bb{\vect k,\bset{\alpha_j}}\cdot D_\nu\bb{\sigma\cdot s+\vect{\sigma}\cdot\vect{k}}
$$
can be written as a finite sum of meromorphic functions and an infinite sum, which is absolutely convergent to a holomorphic function in $P_N=\bset{s\in\C\,:\, \Re(s)> \frac{S_0-\sigma_{N+1}}{\sigma}}$.
\end{lem}
\begin{proof}
Notice that $D_\nu(s)$ is holomorphic in $\mathcal{P}_0 = \{\Re(s)>S_0\}$ and extends meromorphically to $\C$, as it differs from $D(s)$ by the entire function $\sumit n 1 {\nu-1} a_n\cdot c_n^s$.
Define
$$
k_N^*:=\min \bset{k\in\N\,:\,k\ge\frac{\sigma_{N+1}}{\sigma_1}}=\left\lceil{\frac{\sigma_{N+1}}{\sigma_1}}\right\rceil,
$$
and note that if $k_1,\dots,k_N$ satisfy $\sumit j 1 N k_j\ge k_N^*$, then
\[\vect{\sigma}\cdot\vect{k} = \sumit j 1 N \sigma_jk_j \ge \sumit j 1 N \sigma_1k_j\ge \sigma_1k_N^*\ge\sigma_{N+1},\]
and therefore, for $s\in P_N$ (i.e., where $\Re(s)\ge \frac{S_0-\sigma_{N+1}}{\sigma}$), we have 
\[\Re(\sigma\cdot s+\vect{\sigma}\cdot\vect{k})>\sigma\frac{S_0-\sigma_{N+1}}{\sigma}+\sigma_{N+1}=S_0.\]
In other words, \ $P_N\subset\{s\in\C\;:\; \sigma\cdot s+\vect{\sigma}\cdot\vect{k}\in\mathcal{P}_0\}$.
This implies that whenever $\sumit j1 N k_j\ge k_N^*$, the function $s\mapsto D_\nu\bb{\sigma\cdot s+\vect{\sigma}\cdot\vect{k}}$ is holomorphic in $P_N$ and,
$$
D_\nu\bb{\sigma\cdot s+\vect{\sigma}\cdot\vect{k}}=\sumit n \nu \infty a_n\cdot c_n^{\sigma\cdot s+\vect{\sigma}\cdot\vect{k}}.
$$
Next, we shall bound the tail of the sum defining $\widetilde{D}(s)$. Using \eqref{eq:falling}, we have
\begin{align*}
&\sumit k {k_N^*}\infty\abs{\frac{(s)_k}{k!}}\underset{\vect k,\; \sumit j 1 {N} k_j= k}\sum \abs{A\bb{\vect k,\bset{\alpha_j}}}\abs{D_\nu\bb{\sigma\cdot s+\vect{\sigma}\cdot\vect{k}}}\\
&\quad\quad\quad\le\abs s\cdot \sumit k {k_N^*}\infty k^{\abs{s+1}}\underset{\vect k,\; \sumit j 1 {N} k_j= k}\sum \abs{A\bb{\vect k,\bset{\alpha_j}}}\sumit n \nu \infty \abs{a_n\cdot c_n^{\sigma\cdot s+\vect{\sigma}\cdot\vect{k}}}\\
&\quad\quad\quad=\abs s\cdot \sumit n \nu \infty \abs{a_n\cdot c_n^{\sigma\cdot s}}\sumit k {k_N^*}\infty k^{\abs{s+1}}\underset{\vect k,\; \sumit j 1 {N} k_j=k}\sum A\bb{\vect k,\bset{\abs{\alpha_j}}}\cdot |c_n|^{\vect{\sigma}\cdot\vect{k}}.
\end{align*}
To bound the latter, recall that $\delta$, defined in part \ref{item:err} of Lemma \ref{lem:apprx_poly_h}, satisfies \eqref{eq:delta} while $|c_n|<\delta$ for all $n\ge\nu$. We may then write $|c_n|^{\vect{\sigma}\cdot\vect{k}} = \delta^{\vect{\sigma}\cdot\vect{k}}\left(\frac{|c_n|}{\delta}\right)^{\vect{\sigma}\cdot\vect{k}}$, and, since $\frac{\abs{c_n}}\delta<1$ and $\vect{\sigma}\cdot\vect{k}\ge\sigma_{N+1}$ whenever $\sumit j 1 N k_j\ge k_N^*$, we have
\begin{align*}
\underset{\vect k,\; \sumit j 1 {N} k_j=k}\sum A\bb{\vect k,\bset{\abs{\alpha_j}}}|c_n|^{\vect{\sigma}\cdot\vect{k}}&\le \left(\underset{\vect k,\; \sumit j 1 {N} k_j=k}\sum A\bb{\vect k,\bset{\abs{\alpha_j}}}\delta^{\vect{\sigma}\cdot\vect{k}}\right)\left(\frac{\abs{c_n}}{\delta}\right)^{\sigma_{N+1}} \\
&= \left(\sumit j 1 N |\alpha_j|\delta^{\sigma_j}\right)^k\left(\frac{\abs{c_n}}{\delta}\right)^{\sigma_{N+1}} \le 2^{-k}\cdot\delta^{-\sigma_{N+1}}\cdot\abs{c_n}^{\sigma_{N+1}},
\end{align*}
since $\left(\sumit j 1 N |\alpha_j|\delta^{\sigma_j}\right)^k\le 2^{-k}$. We conclude that,
\begin{align*}
&\sumit k {k_N^*}\infty\abs{\frac{(s)_k}{k!}}\underset{\vect k,\; \sumit j 1 {N} k_j= k}\sum \abs{A\bb{\vect k,\bset{\alpha_j}}}\abs{D_\nu\bb{\sigma\cdot s+\vect{\sigma}\cdot\vect{k}}}\\
&\quad\quad\quad=\abs s\cdot \sumit n \nu \infty \abs{a_n\cdot c_n^{\sigma\cdot s}}\sumit k {k_N^*}\infty k^{\abs{s+1}}\underset{\vect k,\; \sumit j 1 {N} k_j=k}\sum A\bb{\vect k,\bset{\abs{\alpha_j}}}|c_n|^{\vect{\sigma}\cdot\vect{k}} \\
&\quad\quad\quad \le \abs s\cdot \sumit n \nu \infty \abs{a_n\cdot c_n^{\sigma\cdot s}}\sumit k {k_N^*}\infty k^{\abs{s+1}}2^{-k}\delta^{-\sigma_{N+1}}|c_n|^{\sigma_{N+1}}\\
&\quad\quad\quad=\frac{\abs s}{\delta^{\sigma_{N+1}}}\cdot \left(\sumit n \nu \infty \abs{a_n\cdot c_n^{\sigma\cdot s+\sigma_{N+1}}}\right)\left(\sumit k {k_N^*}\infty k^{\abs{s+1}}2^{-k}\right).
\end{align*}
The second sum converges independently of $s$, while the first sum converges whenever $\sigma\cdot s+\sigma_{N+1}\in \mathcal{P}_0$, i.e., when $s\in P_N$. We conclude that the sum
$$
\sumit k {k_N^*}\infty \frac{(s)_k}{k!}\underset{\vect k,\; \sumit j 1 {N} k_j= k}\sum A\bb{\vect k,\bset{\alpha_j}}\cdot D_\nu\bb{\sigma\cdot s+\vect{\sigma}\cdot\vect{k}}
$$
is absolutely convergent and converges to a holomorphic function in $P_N$, as a local uniform limit of such functions. Next, the function
$$
\sumit k 0{k_N^*-1}\frac{(s)_k}{k!}\underset{\vect k,\; \sumit j 1 {N} k_j= k}\sum A\bb{\vect k,\bset{\alpha_j}}D_\nu\bb{\sigma\cdot s+\vect{\sigma}\cdot\vect{k}}
$$
is meromorphic as a finite sum of meromorphic functions. Overall, we conclude that
\begin{align*}
\widetilde D(s)&:=\sumit k {0}\infty\frac{(s)_k}{k!}\underset{\vect k,\; \sumit j 1 {N} k_j= k}\sum A\bb{\vect k,\bset{\alpha_j}}D_\nu\bb{\sigma\cdot s+\vect{\sigma}\cdot\vect{k}}\\
&=\sumit k 0{k_N^*-1}\frac{(s)_k}{k!}\underset{\vect k,\; \sumit j 1 {N} k_j= k}\sum A\bb{\vect k,\bset{\alpha_j}}D_\nu\bb{\sigma\cdot s+\vect{\sigma}\cdot\vect{k}}\\
&+\sumit k {k_N^*}\infty\frac{(s)_k}{k!}\underset{\vect k,\; \sumit j 1 {N} k_j= k}\sum A\bb{\vect k,\bset{\alpha_j}}D_\nu\bb{\sigma\cdot s+\vect{\sigma}\cdot\vect{k}}
\end{align*}
is meromorphic in $P_N$ as a finite sum of meromorphic functions and a holomorphic function, thus concluding the proof of the lemma.
\end{proof}

%% file: sections/proof.tex
We begin this section with an observation showing that it suffices to find a sequence of extensions, meromorphic on an increasing sequence of half-planes, in order to generate an extension to the entire complex plane.
\begin{obs}\label{obs:mer_convergece}
Let $\bset{g_m}$ be a sequence of functions satisfying that for every $m$ the function $g_m(z)$ is meromorphic in a half-plane $P_m$, where $P_m\nearrow\C$, while $\left.g_m\right|_{P_1}=\left.g_1\right|_{P_1}$. Then the sequence $\bset{g_m}$ converges locally uniformly to a function, $g$, which is meromorphic in $\C$.
\end{obs}
\begin{proof}
Fix $m>0$ and let $g_k,g_j$ be two functions, defined on $P_m$. By uniqueness of meromorphic functions, if $\left.g_k\right|_{P_1}=\left.g_1\right|_{P_1}=\left.g_j\right|_{P_1}$, then $\left.g_k\right|_{P_m}=\left.g_j\right|_{P_m}$.

For every disk, $B(0,R)\subset\C$, there exists $m$ so that $B(0,R)\subset P_m$. The sequence $\bset{g_j}$ is a sequence of meromorphic functions in $P_m$ and, in particular, in $B(0,R)$, and for every $k\ge m$, we saw that $\left.g_k\right|_{B(0,R)}=\left.g_m\right|_{B(0,R)}$, i.e., the subsequence, $\bset{g_k}_{k=m}^\infty$ is fixed when restricted to $B(0,R)$. We conclude that $\bset{g_k}$ converges locally uniformly to a meromorphic function, $g$, concluding the proof.
\end{proof}

We are finally ready to prove Theorem \ref{thm:extension}. Recall

\begin{thm*}[Theorem \ref{thm:extension}]
Let $D(s):=\sumit n 1 \infty a_n\cdot c_n^s$ be a Dirichlet series with positive bases $c_n$ that converges absolutely in the half-plane $\mathcal{P}_0=\bset{\Re(s)>S_0}$, and extends meromorphically to $\C$. 
Let $\Gamma$ be a bounded sector containing the set $\{c_n\}$, and let $h:\Gamma\rightarrow\C$
be a function satisfying that for every $N$, there exists $C_N>1$ such that 
\[\tag{\ref{eq:cond_h}}
\abs{h(z)-\sumit j1N \alpha_j\cdot  z^{\sigma_j}}\le C_N\cdot \abs z^{\sigma_{N+1}} \text{ for all }z\in \Gamma
\]
where $\{\alpha_j\}\subset\C$ and $\bset{\sigma_j}\subset(0,\infty)$ satisfies $\sigma_j\nearrow \infty$.

Then the sum $\sumit n 1 \infty a_n\cdot \left(c_n^{\sigma}\bb{1+h(c_n)}\right)^s$ extends meromorphically to $\C$ as well.
\end{thm*}
\begin{proof}
Note that, by comparison to the series $\sumit n 1 \infty a_n\cdot c_n^s$, we have that $\sumit n 1 \infty a_n\cdot (c_n^\sigma (1+h(c_n)))^s$ converges absolutely on $P_{-1} = \{s\in\C\;:\;\Re(s)>\frac{S_0}\sigma\}$ and defines a holomorphic function there. 
We will show there exists a sequence of meromorphic functions, $g_N:P_N\rightarrow\C$, that agree with the original sum on $P_{-1}$, and use Observation \ref{obs:mer_convergece}, noting that $P_N\nearrow\C$, to conclude that the sequence converges to a function, $g$, which is meromorphic in $\C$ and agrees with the sum in $P_{-1}$, i.e., a meromorphic extension for the sum in $\C$.

We begin by constructing the sequence $\{g_N\}$ by defining two auxiliary functions. For a fixed $N\in \N$, let $\nu_N$ be a value of $\nu$ defined in Lemma \ref{lem:analytic_tail} (note this depends on $N$), 
and recall that $D_{\nu_N}(s):=\sumit n {\nu_N} \infty a_n\cdot c_n^s$. Define
\begin{align*}
&h_N(s):=\sumit k 1\infty\frac{(s)_k}{k!}\sumit n {\nu_N} \infty a_n\cdot c_n^{\sigma\cdot s}\cdot\mathcal E_k^{(N)}(c_n)\\
&q_N(s)
:=\sumit k {0}\infty\frac{(s)_k}{k!}\underset{\vect k,\; \sumit j 1 {N} k_j= k}\sum A\bb{\vect k,\bset{\alpha_j}}D_{\nu_N}\bb{\sigma\cdot s+\vect{\sigma}\cdot\vect{k}},
\end{align*}
where $\mathcal E_k^{(N)}$ is the function $\mathcal E_k$ defined in Lemma \ref{lem:apprx_poly_h}, which, again, implicitly depends on $N$.
Following Lemma \ref{lem:analytic_tail}, the function $h_N$ is analytic in $P_N$, and following Lemma \ref{lem:mero-comb}, the function $q_N$ is meromorphic in $P_N$, implying that the function $g_N:P_N\rightarrow\C$ defined by 
$$
g_N(s)=\sumit n 1 {\nu_N-1} a_n\cdot c_n^{\sigma\cdot s}\bb{1+h(c_n)}^s+q_N(s)+h_N(s)
$$
is meromorphic in $P_N$ as a finite sum of holomorphic and meromorphic functions in $P_N$.

To conclude the proof, it is left to show that for every $s\in P_{-1}$
\begin{align}\label{eq:new_goal}
g_N(s)=\sumit n 1 \infty a_n\cdot c_n^{\sigma\cdot s}\bb{1+h(c_n)}^s,
\end{align}
i.e., $g_N$ agrees with the original sum on $P_{-1}$. Indeed, note that if $s\in P_{-1}$, then
$$
D_{\nu_N}\bb{\sigma\cdot s+\vect{\sigma}\cdot\vect{k}}=\sumit n {\nu_N} \infty a_n\cdot c_n^{\sigma\cdot s+\vect{\sigma}\cdot\vect{k}}.
$$
Following Lemma \ref{lem:analytic_tail} and Lemma \ref{lem:mero-comb} the sums composing $q_N$ and $h_N$ are absolutely convergent in $P_{-1}$ (excluding a finite number of elements), implying that for $s\in P_{-1}$ we have
\begin{align*}
q_N(s)+h_N(s) &= \sumit k {0}\infty\frac{(s)_k}{k!} \underset{\vect k,\; \sumit j 1 {N} k_j= k}\sum A\bb{\vect k,\bset{\alpha_j}}D_{\nu_N}\bb{\sigma\cdot s+\vect{\sigma}\cdot\vect{k}} + \sumit k {1}\infty\frac{(s)_k}{k!}\sumit n {\nu_N} \infty a_n\cdot c_n^{\sigma\cdot s}\cdot\mathcal E_k(c_n) \\
&=\sumit n {\nu_N} \infty \left(a_n\cdot c_n^{\sigma\cdot s}+\sumit k {1}\infty\frac{(s)_k}{k!}\bb{\underset{\vect k,\; \sumit j 1 {N} k_j= k}\sum A\bb{\vect k,\bset{\alpha_j}} a_n\cdot c_n^{\sigma\cdot s+\vect{\sigma}\cdot\vect{k}}+a_n\cdot c_n^{\sigma\cdot s}\cdot\mathcal E_k(c_n)}\right) \\
&=\sumit n {\nu_N} \infty a_n\cdot c_n^{\sigma\cdot s}\left(1+\underset{\vect k,\; \sumit j 1 {N} k_j= k}\sum A\bb{\vect k,\bset{\alpha_j}}  c_n^{\vect{\sigma}\cdot\vect{k}}+\mathcal E_k(c_n)\right) = \sumit n {\nu_N} \infty  a_n\cdot c_n^{\sigma\cdot s}\bb{1+h(c_n)}^s,
\end{align*}
and hence
\begin{align*}
g_N(s)&=\sumit n 1 {\nu_N-1} a_n\cdot c_n^{\sigma\cdot s}\bb{1+h(c_n)}^s+q_N(s)+h_N(s)\\
&=\sumit n 1 {\nu_N-1} a_n\cdot c_n^{\sigma\cdot s}\bb{1+h(c_n)}^s+\sumit n {\nu_N} \infty  a_n\cdot c_n^{\sigma\cdot s}\bb{1+h(c_n)}^s=\sumit n 1 {\infty} a_n\cdot c_n^{\sigma\cdot s}\bb{1+h(c_n)}^s.
\end{align*}
This implies \eqref{eq:new_goal} and concludes the proof of the theorem.
\end{proof}
\begin{rmk}\label{rmk:entire}
Looking at the proof, we see that if the extension of the original series, $D(s)$, is an entire function, then for every $N$, the function $g_N$ is holomorphic in $P_N$. This implies that the extension of the perturbation, which is the limit of $\bset{g_N}$, is also an entire function, as a local uniform limit of holomorphic functions.
\end{rmk}
\begin{rmk}\label{rmk:poles}
Looking at the proof more closely, in particular at Lemma \ref{lem:mero-comb}, we see that on the half-plane $P_N$, the new meromorphic extension equals to the finite sum of meromorphic functions
\begin{equation}
\label{eq:finite-mero}
\sumit k 0{k_N^*-1}\frac{(s)_k}{k!}\underset{\vect k,\; \sumit j 1 {N} k_j= k}\sum A\bb{\vect k,\bset{\alpha_j}}D_\nu\bb{\sigma \cdot s+\vect{\sigma}\cdot\vect{k}}\quad\bb{\text{where }k_N^*=\left\lceil{\frac{\sigma_{N+1}}{\sigma_1}}\right\rceil}
\end{equation}
up to a function which is holomorphic on $P_N$. In particular, the poles of the new meromorphic extension are contained in the set
$$
\bunion N 1 \infty\bunion k 0 \infty \underset{\bb{k_1,\cdots,k_N}\atop \sumit j1 Nk_j=k}\bigcup \frac{\mathcal P-\sumit j 1 N\sigma_{j-1}\cdot k_j}\sigma,
$$
where $\mathcal P$ is the set of poles of the original meromorphic extension, $D$. 
The containment can be strict; see Remark \ref{rmk:hurwitz}. More precisely, if we know the Laurent series of $D(s)$ up to a holomorphic function (locally) at a pole, $p$, then computing the coefficients of negative powers of $s-p$ arising in the finite sum in \eqref{eq:finite-mero}, one can realise wether $p$ is indeed a pole or not. While this set seems \emph{a priori} incredibly complicated to find, in many cases it is very simple (see, e.g., Examples \ref{thm:geo} and \ref{thm:zeta}).
\end{rmk}

\begin{rmk}
If the function $h$ satisfies additional requirements then the constant $\nu_N$, which depends on both the error term and on the sequence $\bset{c_n}$, could be chosen uniformly in $N$.

One notable example is if $h(z) = \sumit j 1 {N_0}\alpha_j\cdot z^{\sigma_j}$ is a finite linear combination of powers of $z$. In this case, we can take $\nu_N = \nu_{N_0}$ for all $N\ge N_0$, since the constant $\delta$, defined in Lemma \ref{lem:apprx_poly_h}, is independent of $N$ for $N\ge N_0$. Then, $g_N=g_{N_0}$
, and moreover for $\nu=\nu_{N_0}$ we have
$$
g_{N_0}=\sumit n 1 {\nu-1} a_n\cdot c_n^{\sigma\cdot s}\bb{1+h(c_n)}^s+\sumit k {0}\infty\frac{(s)_k}{k!} \underset{\vect k,\; \sumit j 1 {N_0} k_j= k}\sum A\bb{\vect k,\bset{\alpha_j}}\cdot D_{\nu}\bb{\sigma\cdot s+\vect{\sigma}\cdot\vect{k}}.
$$

Another notable example is when $h$ is holomorphic at the origin. In this case, by Taylor's theorem, there exists $\eps_0$ such that the Taylor series, $\sumit n 1\infty \alpha_n\cdot z^n$, converges uniformly and in absolute value in $B(0,\eps_0)$. Following Cauchy's integral formula for the derivative, we have
$$
\abs{\alpha_n}=\abs{\frac{h^{(n)}(0)}{n!}}=\abs{\frac1{2\pi i}\oint_{B(0,\eps_0)}\frac{h(z)}{z^{n+1}}dz}\le C_0\cdot\eps_0^{-n},\quad C_0:=\underset{\abs z=\eps_0}\max\; \abs{h(z)}.
$$
Then for every $z\in B\bb{0,\frac{\eps_0}2}$,
\begin{align*}
\abs{h(z)-\sumit n 1 N a_n\cdot z^n}\le \sumit n{N+1}\infty \abs{a_n}\abs z^n\le  \sumit n{N+1}\infty C_0\cdot\eps_0^{-n}\cdot\abs z^n\le C_0\cdot\eps_0^{-(N+1)}\cdot \abs z^{N+1}\sumit k 0 \infty 2^{-k}=\frac{2C_0}{\eps_0^{N+1}}\cdot \abs z^{N+1}.
\end{align*}
We can therefore take $C$ in the definition of the error term to be $\frac{2C_0}{\eps_0^{N+1}}$.
Noting that, if $\delta<\frac{\eps_0}{2}$,
\begin{align*}
&\sumit j 1 N\abs{\alpha_j}\cdot \delta^j\le \sumit j 1 \infty C_0\cdot \eps_0^{-j}\cdot \delta^j= C_0\cdot  \frac{\frac{\delta}{\eps_0}}{1-\frac{\delta}{\eps_0}} < 2C_0\cdot\frac{\delta}{\eps_0}\quad\quad\text{ and }\\
&C\cdot\delta ^{N+1}\le \frac{2C_0}{\eps_0^{N+1}}\cdot\delta^{N+1}= 2C_0\cdot \bb{\frac\delta{\eps_0}}^{N+1}\le 2C_0\cdot\frac\delta{\eps_0},
\end{align*}
it follows that any $\delta < \eps_0\cdot \min\bset{\frac{1}{4C_0},\frac{1}{2}}$ satisfies $\sumit j 1 N\abs{\alpha_j}\cdot \delta^j < \frac12$ and $C\cdot\delta ^{N+1}<\frac12$.
%
%
As both constants, $C_0$ and $\eps_0$, are independent of $N$,
we conclude that, in this case, $\nu$ can be chosen uniformly as well.

If the sequence $\{\sigma_j\}$ is not a sequence of integers, one can show that, if $\inf_j(\sigma_{j+1}-\sigma_j)>0$, then the condition $C_N\le CM^{\sigma_{N+1}}$ for some $C,M>1$ will lead to the same conclusions. Note that if the sequence $\{\sigma_j\}$ is a subset of $\N$, then this condition is precisely that of being analytic near $0$.
\end{rmk}

%% file: sections/examples.tex
This section deals with two kinds of examples. The first subsection will explore the usefulness of Theorem \ref{thm:extension}. In the second subsection we present a counterexample demonstrating that the asymptotic expansion conditions imposed by Theorem \ref{thm:extension} are necessary. Namely, for every $m$ we can construct a perturbation, $h$, which is smooth on $(0,\eps)$, satisfies $\abs{h(z)}\le |z|^m$on $[0,\eps]$, i.e., condition \eqref{eq:cond_h} holds with $\alpha_j\equiv 0$ and with error $\abs z^m$, but the boundary of the half plane where the sum converges to a holomorphic map forms a natural boundary.

\subsection{On the Usefulness of Theorem \ref{thm:extension} -- Examples:}
\subsubsection{On Almost Geometric Progressions:}
The motivation for this paper was given by Nir Avni who asked us whether $\sumit n 1 \infty\frac1{\bb{2^n+1}^s}$ has a meromorphic extension. This led to the following class of examples:
\begin{thm}\label{thm:geo}
Let $a>1, b\in\C $, and suppose $a^n+b\ne 0$ for all $n\in\mathbb{N}$. Then the Dirichlet series $\sumit n 1 \infty \frac1{\bb{a^n+b}^s}$ has a meromorphic extension to the entire complex plane. Moreover, the poles of the extension are contained in the set
\[\frac{2\pi i}{\log(a)}\mathbb{Z}-\mathbb{N}_0 = \left\{\frac{2\pi i}{\log(a)}j-k\,:\, j,k\in\Z, k\ge 0\right\}.\]
\end{thm}
\begin{proof}
Note that since $a>1$, then
$$
\sumit n 1 \infty\frac1{\bb{a^n}^s}=\sumit n 1 \infty\frac1{\bb{a^s}^n}=\frac 1{a^s-1},
$$
i.e., the Dirichlet series $\sumit n 1 \infty\frac1{\bb{a^n}^s}$ has a meromorphic extension to the entire complex plane. Next
$$
\frac1{a^n+b}=\frac1{a^n}\bb{\frac1{1+\frac b{a^n}}}=\frac1{a^n}\bb{1-\frac{\frac b{a^n}}{1+\frac b{a^n}}}=c_n\bb{1+h(c_n)}
$$
for $c_n:=\frac1{a^n}$ and $h(z)=-\frac{b\cdot z}{1+b\cdot z}$, which is holomorphic in a neighbourhood of the origin, which depends on $\abs b$. Following Theorem \ref{thm:extension}, the Dirichlet series $\sumit n 1 \infty \frac1{\bb{a^n+b}^s}$ has a meromorphic extension to the entire complex plane. Finally, following Remark \ref{rmk:poles}, since $h(z)$ admits an asymptotic expansion where the $\sigma_j$ are positive integers, we see that all poles of the extension must be of the form $p-k$ where $p$ is a pole of $\sumit n 1 \infty \frac{1}{(a^n)^s} = \frac{1}{(a^s-1)}$ and $k$ is a nonnegative integer, and since the poles of $\frac{1}{a^s-1}$ are located at $s\in\frac{2\pi i}{\log(a)}\Z$, the  poles of the extension must be in $\frac{2\pi i}{\log(a)}\Z-\N_0$, as claimed.
\end{proof}

A different kind of almost arithmetic progressions is when $c_n$ is a sum of several arithmetic progressions:
\begin{thm}\label{thm:examp_23}
Let $a_1,a_2,\cdots,a_d\in\bb{0,1},\;\alpha_1,\cdots,\alpha_d>0$. The Dirichlet series $\sumit n 1 \infty \bb{\sumit j 1 d \alpha_j\cdot a_j^{n}}^s$ has a meromorphic extension to the entire complex plane.
\end{thm}
In particular, Dirichlet series such as $\sumit n 1 \infty \bb{2^{-n}+3^{-n}}^s$ have meromorphic extensions to the entire complex plane.
\begin{proof}
Without loss of generality $a_1>a_2>\cdots>a_d$ and $\alpha_1=1$. Since for every $a>0$, 
$$
a^{-n}=\exp\bb{-n\log (a)}=\exp\bb{-n\log (a)\cdot\frac{\log(b)}{\log(b)}}=\bb{b^{-n}}^{\frac{\log(a)}{\log(b)}},
$$
we can rewrite $\widehat{c}_n$ as
$$
\widehat{c}_n= a_1^{n}\bb{1+\sumit j 2 d \alpha_j\cdot \bb{a_1^{n}}^{\frac{\log(a_j)}{\log(a_1)}}}= a_1^{n}\bb{1+h(a_1^n)}
$$
for $h(z)=\sumit j 2 d \alpha_j\cdot z^{\sigma_j}$ where $\sigma_j:={\frac{\log(a_j)}{\log(a_1)}}$, is a monotone increasing sequence. Applying Theorem \ref{thm:extension} we conclude that since $\sumit n 1 \infty  (a_1^n)^s
=\frac{a_1^s}{1-a_1^s}$ has a meromorphic extension to the entire complex plane, so does $\sumit n 1 \infty \bb{\sumit j 1 d \alpha_j\cdot a_j^{n}}^s$.
\end{proof}
\subsubsection{An Almost Zeta Function}
\begin{thm}\label{thm:zeta}
Let $\alpha,\beta\in\R$ satisfy $n+\alpha\cdot n^\beta\ne 0$ for all $n\in\N$, as well as $\alpha>0$ if $\beta\ge 1$. Then the 
Dirichlet series
$$
\sumit n 1 \infty\frac1{\bb{n+\alpha\cdot n^\beta}^s}
$$
has a meromorphic extension to the entire complex plane. Moreover, the set of poles of this series, $\mathcal P$ satisfies
$$
\mathcal P
	\begin{cases}
		=\bset{1}&, \beta=1\\
		\subseteq 1-(1-\beta)\N_0&, \beta<1\\
		\subseteq \frac1\beta-\frac{1-\frac1\beta}\beta\N_0&, \beta>1
	\end{cases}.
$$
\end{thm}

\begin{proof}
Recall that
$$
\sumit n 1 \infty\frac1{n^s}=\zeta(s),
$$
where $\zeta(s)$ extends meromorphically to $\C$ with a single pole at $s=1$. Note that if 
$\beta=1$ then 
$\sumit n 1 \infty\frac1{\bb{n+\alpha\cdot n^\beta}^s}=(1+\alpha)^{-s}\sumit n 1 \infty \frac1{n^s}=(1+\alpha)^{-s}\cdot\zeta(s)$. 

Next, if $\beta<1$ then
$$
\frac1{n+\alpha\cdot n^\beta}=\frac1{n}\cdot\frac1{1+\alpha\cdot n^{-\bb{1-\beta}}}=\frac1{n}\cdot\bb{1-\frac{\alpha\cdot n^{-\bb{1-\beta}}}{1+\alpha\cdot n^{-\bb{1-\beta}}}}=c_n\bb{1+h(c_n)}
$$
if $c_n=\frac1n$ and $h(z)=-\frac{\alpha\cdot z^{1-\beta}}{1+\alpha\cdot z^{1-\beta}}$. 
Noting that
$h(z) = \sumit j 1 \infty (-\alpha)^j z^{(1-\beta)j}$ with the sum converging absolutely on the region $\alpha\cdot |z|^{1-\beta}<1$ for $\beta<1$, we see that $\sumit n 1 \infty\frac1{\bb{n+\alpha\cdot n^\beta}^s} = \sumit n 1 \infty c_n^s(1+h(c_n))^s$
has a meromorphic extension to the entire complex plane. Since $h$ admits an asymptotic expansion where the $\sigma_j$ are integer multiplies of $(1-\beta)$, by Remark \ref{rmk:poles}, the poles of the extension are contained in $1-(1-\beta)\N_0$.

Finally, if $\beta>1$ then 
$$
\frac1{n+\alpha\cdot n^\beta}=\frac1{\alpha \cdot n^\beta}\cdot\frac1{1+\alpha^{-1}\cdot n^{-\bb{\beta-1}}}=\frac1{\alpha\cdot n^\beta}\cdot\bb{1-\frac{\alpha^{-1}n^{-\bb{\beta-1}}}{1+\alpha^{-1}\cdot n^{-\bb{\beta-1}}}}=\alpha^{-1}\cdot c_n\bb{1+h(c_n)}
$$
if $c_n=\frac1{n^\beta}$ and $h(z)=-\frac{\alpha^{-1}\cdot z^{1-\frac1\beta}}{1+\alpha^{-1}\cdot z^{1-\frac1\beta}}$. 
Note that
$$
\sumit n 1 \infty\frac1{\bb{n^\beta}^s}=\sumit n 1 \infty\frac1{n^{\beta\cdot s}}=\zeta\bb{\beta\cdot s}
$$
extends to a meromorphic function with a single pole at $s=\frac1\beta$. We conclude that since $h(z) = \sumit j 1 \infty (-\alpha)^{-j}z^{(1-\frac1\beta)j}$ with the sum converging absolutely on $\frac{|z|^{1-\frac1\beta}}{\alpha}<1$ for $\beta>1$, then $\sumit n 1 \infty\frac{1}{(n+\alpha\cdot n^\beta)^s} = \alpha^{-s}\sumit n 1 \infty c_n^s(1+h(c_n))^s$ has a meromorphic extension to the entire complex plane. Since $h$ admits an asymptotic expansion where the $\sigma_j$ are integer multiplies of $\bb{1-\frac1\beta}$, by Remark \ref{rmk:poles}, the poles of the extension are contained in $\frac1\beta-\frac{1-\frac1\beta}\beta\N_0$.
\end{proof}
\begin{rmk}
\label{rmk:hurwitz}
If $\beta=0$, then
\[\sumit n 1 \infty \frac{1}{\bb{n+\alpha\cdot n^0}^s} = \sumit n 0 \infty \frac{1}{\bb{n+1+\alpha}^s} = \zeta(s,1+\alpha),\]
where $\zeta(s,a)$ is the \emph{Hurwitz zeta function}. This example thus recovers the fact that $\zeta(s,a)$ extends to a meromorphic function, with potential poles at $1-\N_0$. However, a careful analysis of the terms appearing in \eqref{eq:finite-mero} shows that the poles that could appear at $-\N_0$ are in fact cancelled. That is, if we write $\zeta(s,a) = \sumit k 0 N p_k(s)\zeta(s+k)$ modulo a holomorphic function on $P_N$ for some polynomials $p_k(s)$, then $p_k(1-k)=0$ for all $k\ge 1$, thus cancelling the simple pole of $\zeta(s+k)$ at $s=1-k$. This recovers the fact that $\zeta(s,a)$ in fact only has a simple pole at $s=1$.
\end{rmk}

A different class of examples arises from looking at sums $\sum \frac1{p(n)^s}$ for a generic polynomial $p$.
\begin{thm}
Let $p(z)$ be a monic nonzero degree $d$ polynomial, and define
$$
\Sigma(s):=\sumit n {N} \infty\frac1{\bb{p(n)}^s}
$$
where $N$ is sufficiently large so that $p(n)\ne 0$ for all $n\ge N$. Then the Dirichlet series, $\Sigma$, has a meromorphic extension to the entire complex plane, and moreover its poles are contained in $\frac{1}{d}\left(1-\N_0\right).$
\end{thm}
\begin{proof}
Writing $p(z)=z^d+\sumit j 0 {d-1}a_jz^j$, we proceed similarly as in the previous example and write
$$
\frac1{p(n)}=\frac1{n^{d}\bb{1+\sumit j 0 {d-1}a_j\cdot n^{j-d}}}=\left(\frac{1}{n}\right)^d\bb{1-\frac{\sumit j 0 {d-1}a_j\cdot \bb{\frac1{n}}^{d-j}}{1+\sumit j 0 {d-1}a_j\cdot \bb{\frac1{n}}^{d- j}}}.
$$
As before we may apply Theorem \ref{thm:extension} and Remark \ref{rmk:poles} with $c_n=\frac1{n}$, $\sigma=d$, and the function $h(z):=-\frac{\sumit j 0 {d-1}a_j\cdot z^{d-j}}{1+\sumit j 0 {d-1}a_j\cdot z^{d-j}}$, which is holomorphic in $B(0,\eps_p)$ for some $\eps_p>0$.
\end{proof}

\subsubsection{Dirichlet L-series}
A notable example where Remark \ref{rmk:entire} might be useful is for perturbations of \emph{Dirichlet L-series}. These are sums of the form $\sumit n 1\infty\frac{\chi(n)}{n^s}$, where $\chi(n)$ are \emph{Dirichlet characters  of modulus $m$} (for some $m\in\N$), namely functions $\chi:\N\to\C$ satisfying the following properties for every $a,b\in\N$:
\begin{enumerate}
\item completely multiplicative: $\chi(a\cdot b)=\chi(a)\cdot \chi(b)$.
\item $\chi(a)=0\iff gcd(m,a)>1$.
\item periodic with period $m$: $\chi(a+m)=\chi(a)$.
\end{enumerate}
In our setup this translates to $c_n=\frac1n$ and $a_n=\chi(n)$. If $\chi(n)=1$ for all $n$, this reduces to the \emph{Riemann zeta function}, which extends meromorphically to $\C$ with a single pole at $s=1$. 
For \emph{non-principal} $\chi$, 
the series extends to an \emph{entire} function on $\C$, allowing us to apply Remark \ref{rmk:entire} and obtain an entire extension for perturbations of such series.
For more information see \cite[chapter 16]{IR}.

\subsubsection{A Rough Perturbation - Still Good Enough:}\label{subsec:rough}
Let $W(t)$ denote the Weierstrass function
$$
W(t):=\sumit n 0\infty \bb{\frac56}^n\cos\bb{7^n\pi\cdot t}.
$$
The function $W$ is nowhere differentiable and bounded by 6. Next, we define the function $h(t):=e^{-\frac1{t^2}}\cdot W(t)$. Then for every $N$ there exists $C_N$ so that for every $\abs t<1$,
$$
\abs{h(t)}<C_N\abs t^N.
$$
If the sum $\sumit n 1 \infty a_n\cdot c_n^s$ has a meromorphic extension to the entire complex plane, then the perturbed sum $\sumit n 1 \infty a_n \cdot c_n^s\bb{1+h(c_n)}^s$ satisfies the requirements of Theorem \ref{thm:extension}, and therefore it also has a meromorphic extension to the entire complex plane, despite the fact that $h$ is far from being differentiable at any point but the origin. Moreover, since $h$ has an asymptotic expansion which is identically zero, it follows that the difference between the perturbed series and the original series, $\sumit n 1 \infty a_n\cdot c_n^s\bb{1+h(c_n)}^s - \sumit n 1 \infty a_n\cdot c_n^s$, extends to an \emph{entire} function on $\C$. To see this, we recall the definition of the function $q_N$ in the proof of Theorem \ref{thm:extension},
$$
q_N(z)=\sumit k {0}\infty\frac{(s)_k}{k!}\underset{\vect k,\; \sumit j 1 {N} k_j= k}\sum A\bb{\vect k,\bset{\alpha_j}}D_{\nu_N}\bb{\sigma\cdot s+\vect{\sigma}\cdot\vect{k}}=D_{\nu_N}(\sigma\cdot s),
$$
 since $\alpha_j=0$ for all $j$ and therefore $A\bb{\vect k, \bset{\alpha_j}}=0$ for all $\vect{k}$ except $\vect{k} = (0,\dots,0)$. This implies that for $s\in P_N$
\begin{align*}
\sumit n 1 \infty a_n \cdot c_n^s\bb{1+h(c_n)}^s - \sumit n 1 \infty &a_n\cdot c_n^s= g_N(s) - \sumit n 1 \infty a_n\cdot c_n^s\\
&=\sumit n 1 {\nu_N-1} a_n\cdot c_n^s\bb{1+h(c_n)}^s+q_N(s)+h_N(s)- \sumit n 1 {\nu_N-1} a_n\cdot c_n^s - \sumit n {\nu_N} \infty a_n\cdot c_n^s\\
&=\sumit n 1 {\nu_N-1} a_n\cdot c_n^s\left( \bb{1+h(c_n)}^s-1\right)+D_{\nu_N}(s)+h_N(s) - \sumit n {\nu_N} \infty a_n\cdot c_n^s \\
&= \sumit n 1 {\nu_N-1} a_n\cdot c_n^s\left( \bb{1+h(c_n)}^s-1\right)+h_N(s),
\end{align*}
i.e., the difference is a sum of $h_N$ and a finite sum of holomorphic functions, and hence is itself holomorphic in the half-plane $P_N$, as $h_N$ is holomorphic in $P_N$. We conclude that the difference is an entire function, as a local uniform limit of entire functions.

\subsection{On the Sharpness of Theorem \ref{thm:extension} -- The Proof of Theorem \ref{thm:example}:}
In this subsection, given any $m \in \mathbb{N}$, we construct a function $h$ that is smooth on $(0, \infty)$ and bounded from above by $\abs z^m$, for which the Dirichlet series 
$\sumit n 1 \infty \left(e^{-n}(1 + h(e^{-n}))\right)^s$ 
converges in the half-plane $\{\operatorname{Re}(s) > -m\}$, but exhibits a natural boundary along the line $\{\operatorname{Re}(s) = -m\}$. Note that this series is a perturbation of the standard Dirichlet series 
$\sumit n 1 \infty e^{-ns} = \frac{1}{e^s - 1}$ which, in contrast to our perturbed series, converges
to a holomorphic function in $\{\operatorname{Re}(s) > 0\}$ and admits a meromorphic extension to the entire complex plane.
The contrast between the behaviours of the two series is due to the function $h$ not satisfying the asymptotic expansion requirement \eqref{eq:cond_h} to all orders.

Recall from the introduction that there exists a sequence $\{a_n\}$, bounded in absolute value by $1$, such that the power series
$\sumit n 0 \infty a_n\cdot z^n$
has a natural boundary at $\bset{|z|=1}$. We use this sequence, $\{a_n\}$, to define our counterexample. 
Let $\chi:\R\to[0,1]$ denote the bump function,
$$
\chi(t)=	\begin{cases}
			\exp\bb{\frac{t^2}{4t^2-1}}&, \abs t<\frac12\\
			0&, \text{otherwise}
		\end{cases}.
$$
It is a smooth function supported in $\bb{-\frac12,\frac12}$ satisfying $\chi(0)=1$ and $0\le\chi(t)\le 1$ for all $t$. Define the function $h(z):(0,\infty)\rightarrow\C$ by
$$
h(z):=z^m\cdot \sumit n 1 \infty a_n\cdot\chi(e^n\cdot z-1).
$$
Note that since $\chi$ is supported on $\bb{-\frac12,\frac12}$, for every $z$ there exists at most one $n$ such that $a_n\cdot\chi(e^n\cdot z-1)\neq 0$. Consequently, $h(z)$ is smooth on $(0,\infty)$, $\abs{h}\le 1$ and moreover, for every $n$, $h(e^{-n})=a_n\cdot e^{-nm}$.

Next, as we saw before, using the Taylor expansion of $z\mapsto\bb{1+z}^s$ centered at $z=0$, we see that
\[(1+z)^s = 1+sz + R_1(z;s),\]
where, if $\abs z<1$, then $s\to R_1(z;s) = (1+z)^s-(1+sz)$ is an entire function in $s$, satisfying $|R_1(z;s)|\le C(s)z^2$ for some continuous function $s\mapsto C(s)$ (locally uniformly in $z$). Then
\begin{align*}
    \sumit n 1 \infty c_n^s(1+h(c_n))^s = \sumit n 1 \infty e^{-ns}(1+h(e^{-n}))^s &= \sumit n 1 \infty e^{-ns}(1+s\cdot h(e^{-n})+R_1(h(e^{-n});s)) \\
    &= \sumit n 1 \infty e^{-ns} + s\sumit n 1 \infty e^{-ns}h(e^{-n}) +\sumit n 1 \infty e^{-ns}R_1(h(e^{-n});s),
\end{align*}
since each summand converges in absolute value in $\bset{\Re(s)>-m}$. The first term is $\sumit n 1 \infty e^{-ns} = \frac{1}{e^s-1}$, which admits a meromorphic extension to $\mathbb{C}$. The elements of the last summand satisfy
\[|e^{-ns}R_1(h\bb{e^{-n}};s)|\le C(s)e^{-ns}|h(e^{-n})|^2 \le C(s)e^{-ns}|e^{-n}|^{2m} = C(s)e^{-n(s+2m)},\]
implying that the series converges absolutely in $\bset{\Re(s)>-2m}$ and, in particular, the sum extends holomorphically to $\bset{\Re(s)>-2m}$, since $s\mapsto e^{-ns}R_1\bb{h(e^{-n});s}$ is entire for every $n$.

Finally, to explore the middle term, we use $h(e^{-n}) = a_n\cdot e^{-nm}$ to write
\[
    s\sumit n 1 \infty e^{-ns}h(e^{-n}) = s\sumit n 1 \infty a_ne^{-n(s+m)}.
\]
We see that this sum converges in $\bset{\Re(s)>-m}$. However, the sum does not admit an extension on any larger open set, since if it did, setting $z = e^{-(s+m)}$, we would generate an extension for the power series $\sumit n 0 \infty a_n\cdot z^n=a_0+\sumit n 1 \infty a_n\cdot z^n$ beyond the unit circle, which is a contradiction to Hadamard's gap theorem.

Thus,
\[\sumit n 1 \infty c_n^s(1+h(c_n))^s = \sumit n 1 \infty e^{-ns} + s\sumit n 1 \infty e^{-ns}h(e^{-n}) + \sumit n 1 \infty e^{-ns}R_1(h(e^{-n});s).\]
is the sum of two terms which {\bf do} admit meromorphic extensions to $\{\Re(s)>-2m\}$ and a term which admits {\bf no} extensions past $\{\Re(s)>-m\}$, and hence the perturbed sum cannot admit an extension past $\{\Re(s)>-m\}$; in particular, it does not admit a meromorphic extension to $\mathbb{C}$.

\begin{rmk}
The very first counterexample we constructed (for ourselves) used the function $h(z)=-\frac1{4\log(z)}$. Using the polylogarithm properties, in particular, the fact that if $\abs s<2\pi$ then
$$
Li_k(e^{-s})=\frac{(-1)^{k-1}\cdot s^{k-1}}{(k-1)!}\bb{H_{k-1}-\log(s)}+\underset{n=0\atop n\neq k-1}{\overset\infty\sum}\frac{\zeta(k-n)}{n!}(-1)^n s^n,\quad\quad\text{ where } H_0=0,\quad H_n=\sumit { j} 1 n\frac1{ j},
$$
we were able to show that the sum  $\sumit n 1 \infty e^{-ns}\bb{1+h(e^{-n})}^s$ decomposes into three sums: two sums that can be extended meromorphically  in a neighbourhood of the origin, and another sum that decays logarithmically near the origin, and therefore cannot be extended meromorphically. We concluded that the Dirichlet sum $\sumit n 1 \infty e^{-ns}\bb{1+h(e^{-n})}^s$ cannot be extended meromorphically to a neighborhood of the origin in $\C$. We then used the fact that any irrational rotation has a dense orbit to show that the function $h(z)=-\frac{e^{-i\alpha\cdot \log(z)}}{4\log(z)}$ for some $\alpha\nin\pi\Q$ gives rise to a Dirichlet sum that has a natural boundary along the line $\bset{\Re(s)=0}$. We did not include the estimates proving this rigorously for sake of brevity; however, we mention the ideas here as they might be of independent interest.
\end{rmk}